%Paper: funct-an/9211002
%From: arveson@math.berkeley.edu (Bill Arveson)
%Date: Sun, 22 Nov 92 14:59:38 PST
%Date (revised): Tue, 24 Nov 92 07:18:54 PST
%Date (revised): Sun, 17 Jan 93 17:13:04 PST

%
%NOTE: This document requires AMSTeX 2.0 or 2.1, the
%corresponding version of the amsppt style macros,
%and the font package amsfonts 2.0 or 2.1. If you try
%to typeset it using an older version of TeX, it
%probably won't work.  You con probably get by
%without the fonts, but you will surely need
%the font metrics that go with amsfonts 2.0 or 2.1.
%
%It typesets without errors on my system.
%
%NOTE: The entire AMSTeX 2.1 package can be obtained
%(free) from the AMS's ftp site in Providence. If you
%are running AMSTeX 2.0, the AMS strongly recommends that
%you upgrade to 2.1.  In particular, this will resolve some
%difficulties that 2.0 has with some fonts, and will
%improve backward compatibility with older versions of TeX.
%
%
%
%This paper was formerly titled "Matrix approximations
%of operator spectra".  The only nontrivial change
%(aside from the new title) is a new abstract.
%
%MACINTOSH USERS:
%I've written a Mac program which illustrates the main
%results of this paper.  You need a Mac II, IIx, IIcx, IIci
%IIfx, Quadra or an SE/30 to use it.  If you want a copy of the
%program + documentation, I can send it binhexed via email.
%(if you don't know what I'm talking about then you
%probably won't be interested in that option.)
%Otherwise, send me a blank Mac floppy for a "hard" copy.
%
%%%%%%%%%%%%%%%%%%%  BEGINNING OF NumLinAlg.tex  %%%%%%%%%%%
%
\input amstex

%MACROS
\def\cstar{$C^*$-algebra}
\def\trace{\roman{trace}}
\def\I{\bold1}
\def\norm{\parallel}
\def\deg{\roman{deg }}
%
%AMSTeX HEADERS
\documentstyle{amsppt}
\loadbold
\magnification=\magstep 1
\topmatter
\title
$C^*$-algebras and numerical linear algebra
\endtitle
\author William Arveson
\endauthor
\affil
Department of Mathematics\\
University of California\\
Berkeley, CA 94720 USA
\endaffil
\email
arveson\@math.berkeley.edu
\endemail
%
%\date
%11 February 1992
%\enddate
%
%
\thanks
This research was supported in part by
NSF grant DMS89-12362
\endthanks
\abstract
Given a self adjoint operator $A$ on a Hilbert space,
suppose that that one wishes to compute the spectrum of
$A$ numerically.  In practice, these problems
often arise in such a way that the matrix of $A$ relative
to a natural basis is ``sparse".
For example, doubly infinite tridiagonal matrices are
usually associated with discretized second order differential operators.
In these cases it is easy and natural to compute the eigenvalues of large
$n\times n$ submatrices of the infinite operator matrix, and
to hope that if $n$ is large enough then the resulting distribution
of eigenvalues will give a good approximation to
the spectrum of $A$.

While this hope is often realized in practice it
often fails as well, and it can fail in spectacular ways.  The
sequence of eigenvalue distributions may not converge as $n\to\infty$,
or they may converge to something that has little to do
with the original operator $A$.  At another level, even the
meaning of `convergence' has not been made precise in general.
In this paper we determine the proper general
setting in which one can expect convergence, and
we describe the asymptotic behavior of the
$n\times n$ eigenvalue distrubutions in all but the most
pathological cases.  Under appropriate hypotheses we establish
a precise limit theorem which shows how the spectrum of
 $A$ is recovered from the sequence of eigenvalues of the
$n\times n$ compressions.

In broader terms, our results have led us to the conclusion
that {\sl numerical problems involving infinite
dimensional operators require a reformulation in terms of}
\cstar s.  Indeed, it is only when the single operator
$A$ is viewed as an element
of an appropriate \cstar\ $\Cal A$ that one can see the precise
nature of the limit of the $n\times n$ eigenvalue distributions; the
limit is associated with a tracial state on $\Cal A$.
Normally, $\Cal A$ is highly
noncommutative, and in the main applications it is a simple
\cstar\ having a unique tracial state.
We obtain precise asymptotic results for the case where $A$
is the discretized Hamiltonian of a one-dimensional quantum system
with arbitrary potential.
\endabstract
\endtopmatter
\vfill\eject
%Replace \pagebreak below with the line above
%to fill the lower part of the title page with
%space, rather than stretching it.
\pagebreak
\document
\subheading{1. Introduction}

There are efficient algorithms available for computing the spectrum of
self-adjoint $n\times n$ matrices.  Using a good desktop computer, one can
find all of the eigenvalues of a $100\times 100$ Toeplitz matrix in a few
seconds.  It is natural to ask how these finite-dimensional methods might
be applied to compute the spectrum of the discretized Hamiltonian of a
one dimensional quantum system, or
a self-adjoint operator in an irrational rotation \cstar.
Considering the burgeoning capabilities of modern computers,
it is unfortunate
that there are apparently no guidelines available which indicate how
one might proceed with such a project except in a few isolated cases.
For example, we were unable to determine from the literature if such a
program is even {\it possible} in general.

The purpose of this paper is to develop some general
methods for computing the spectrum (more precisely, the essential spectrum) of
a self-adjoint operator in terms of the eigenvalues of a sequence of
finite dimensional matrix approximations.  We will see that there is
a natural notion of {\it degree} of an operator (relative to a
filtration of the underlying Hilbert space into finite dimensional
subspaces) which generalizes the
classical idea of band-limited matrices.  More significantly, we will
find that this kind of numerical analysis leads one to a broader
context, in which the operator is seen as an element of a \cstar\
which is frequently simple and possesses at least one tracial state.
The limit distributions arising from the finite dimensional
eigenvalue distributions correspond to traces on this \cstar.

Our intention is not to discuss
algorithms {\it per se}, but rather to present an effective program for
utilizing the finite dimensional techniques to compute infinite dimensional
spectra, along with a precise description of a broad class of operators
for which it will succeed.  Thus, we deal with issues arising from
operator theory and operator algebras; we do not address questions
relating directly to numerical analysis.

In more concrete terms, given an orthonormal basis $\{e_1, e_2, \dots\}$
for a Hilbert space $H$ and an operator $A\in \Cal B(H)$, we may
consider the associated infinite matrix $(a_{ij})$
$$
a_{ij} = <Ae_j,e_i>,\qquad i,j = 1, 2, \dots,
$$
and one can attempt to approximate the spectrum of $A$ by calculating the
eigenvalues of large $n\times n$ submatrices
$$
A_n =
\pmatrix
a_{11}&a_{12}&\hdots&a_{1n}\\
a_{21}&a_{22}&\hdots&a_{2n}\\
\vdots&\vdots&\ddots&\vdots\\
a_{n1}&a_{n2}&\hdots&a_{nn}
\endpmatrix.
$$
This program fails dramatically in general.  For example, if $A$ is the
simple unilateral shift (defined on such a basis by $A e_n = e_{n+1}, n
\geq 1$), then the spectrum of $A$ is the closed unit disk and the essential
spectrum of $A$ is the unit circle, while $\sigma(A_n) = \{0\}$
for every $n$ because each $A_n$ is a nilpotent matrix.

Examples like this have caused operator theorists to view matrix approximations
with suspicion, and to routinely look elsewhere for effective means of
computing spectra.  On the other hand, since the 1920s
there has been evidence that for {\it self-adjoint} operators,
this program can be successfully carried out.

Consider for example the following theorem of Szeg\"o.  Let $f$ be a {\it real}
function in $L^\infty[-\pi,\pi]$ and let $A$ be the multiplication operator
$$
A\xi(x) = f(x)\xi(x),\qquad \xi \in L^2[-\pi,\pi].
$$
With respect to the familiar orthonormal basis $\{e_n : n\in\Bbb Z\}$,
$e_n(x) = e^{inx}$, the matrix of $A$ is a doubly-infinite Laurent
matrix whose principal $n\times n$ submatrices are finite self-adjoint
Toeplitz matrices
$$
A_n =
\pmatrix
a_0&a_{-1}&\hdots&a_{-n+1}\\
a_1&a_0&\hdots&a_{-n+2}\\
\vdots&&\ddots&\vdots\\
a_{n-1}&\hdots&&a_0
\endpmatrix
,
$$
the numbers $a_k$ being the Fourier coefficients of $f$.  Let
$\{\lambda_1, \lambda_2,\dots, \lambda_n\}$ be the eigenvalues of $A_n$,
repeated according to multiplicity.  Szeg\"o's
theorem asserts that for every continuous function $u:\Bbb R\to\Bbb R$
we have
$$
\lim_{n\to\infty}{1 \over n}(u(\lambda_1) + u(\lambda_2) +\dots+ u(\lambda_n))
	   = {1 \over 2\pi}\int_{-\pi}^{\pi} u(f(x))\, dx .\leqno{1.1}
$$
See \cite{8}\cite{10}; a discussion of Szeg\"o's theorem can be found in
\cite{11},
and \cite{6} contains a meticulous discussion of Toeplitz operators.

Let us reformulate Szeg\"o's theorem in terms of weak$^*$-convergence
of measures.  Let $a$ and $b$ be the essential
$\roman{inf}$ and essential $\roman{sup}$ of $f$
$$
\align
a &= \operatornamewithlimits{ess\ inf}_{-\pi\leq x\leq\pi}f(x) \\
b &= \operatornamewithlimits{ess\ sup}_{-\pi\leq x\leq\pi}f(x),
\endalign
$$
and let $\mu$ be the probability measure defined on $[a,b]$ by
$$
\mu(S) = {1\over{2\pi}}m\{x: f(x)\in S\},
$$
$m$ denoting Lebesgue measure on $[-\pi,\pi]$.  Then $1.1$ asserts
that the sequence of discrete probability measures
$$
\mu_n = {1\over n}(\delta_{\lambda_1} + \delta_{\lambda_2} + \dots
+\delta_{\lambda_n}),
$$
$\delta_\lambda$ denoting the unit point mass at $\lambda$, converges to
$\mu$ in the weak$^*$-topology of the dual of $C[a,b]$,
$$
\lim_{n\to\infty}\int_a^b u(x)\,d\mu_n(x) = \int_a^b u(x)\,d\mu(x),\qquad u\in
C[a,b].
$$
The spectrum of $A$ coincides with its essential spectrum in this case, and is
exactly the closed support of the measure $\mu$.

Once one knows that the
measures $\mu_n$ converge in this way to a measure supported
exactly on $\sigma_e(A)$, one may draw
rather precise conclusions about the rate at which the eigenvalues
of $A_n$ accumulate at points of the essential spectrum.
For example, suppose that we choose a point $\lambda\in\sigma_e(A)$.
Then for every open interval $I$ containing $\lambda$ we have $\mu(I)>0$.
If we choose positive numbers $\alpha$ and $\beta$ which are close to
$\mu(I)$ and satisfy $\alpha < \mu(I) < \beta$,
then in the generic case where $A$ has no
point spectrum we may conclude that
$$
n\alpha < \#(\sigma(A_n) \cap I) < n\beta
$$
when $n$ is sufficiently large.  Of course, the symbol $\#$ means the
number of eigenvalues in the indicated set, counting repeated eigenvalues
according to the multiplicity of their occurrence.
In a similar way, the density of eigenvalues of $A_n$ tends to become
increasingly sparse around points
of the complement of the essential spectrum.

In this paper we will show that, while convergence of this type cannot be
expected in general, it does persist in a much broader context than
the above.  Actually, this work was initiated in order to develop
methods for computing the essential spectrum of certain tridiagonal operators
which are associated with the discretized Hamiltonians of
one-dimensional quantum mechanical systems (see section 5).
These operators serve to illustrate the
flavor of our more general results.  Starting with a real valued continuous
function $V$ of a real variable, we fix a real number $\theta$ and
consider the tridiagonal operator $T$ defined in terms of a
bilateral orthonormal basis $\{e_n: n\in \Bbb Z\}$ by
$$
Te_n = e_{n-1} + d_ne_n + e_{n+1},
$$
where the diagonal sequence is defined by $d_n = V(\sin(n\theta))$, $n\in \Bbb
Z$.
The compression $T_n$ of $T$ to the span of $\{e_k: -n\leq k\leq n\}$
is a $(2n+1)\times(2n+1)$ tridiagonal matrix whose eigenvalues can be computed
very rapidly using modern algorithms.  For example, we show that if $\mu_n$ is
the discrete probability measure defined by the eigenvalue distribution of
$T_n$,
$$
\mu_n =
{1\over{2n+1}}(\delta_{\lambda_1}+\delta_{\lambda_2}+\dots\delta_{\lambda_{2n+1}}),
$$
then there is a probability measure $\mu$ {\it having closed support}
$\sigma_e(T)$,
such that for every continuous function $f:\Bbb R\to\Bbb R$ we have
$$
\lim_{n\to\infty}\int f(x)\,d\mu_n(x) = \int f(x)\,d\mu(x).
$$
See formula 5.1.

In these cases the limiting measure $\mu$ can be identified in rather concrete
terms.  Briefly, if $\theta/\pi$ is irrational then the pair of unitary
operators
$$
\align
Ue_n&=e_{n+1}\cr
Ve_n&=e^{in\theta}e_n
\endalign
$$
generate an irrational rotation \cstar\ $\Cal A_\theta$.  $\Cal A_\theta$
has an essentially unique representation
in which its weak closure is the $II_1$ factor $R$.  $T$ becomes
a self-adjoint operator in $R$, and the measure $\mu$ is defined by
$$
\mu(S) = \tau(E_T(S)),\qquad S\subseteq\Bbb R
$$
where $\tau$ is the normalized trace on $R$ and $E_T$ is the spectral measure
of $T$.

In sections 2 and 3 we consider a rather general question.  Given a
self-adjoint
operator $A$ and an orthonormal basis for the underlying Hilbert space,
let $(a_{ij})$ be the matrix of $A$ relative to this basis.  How must one
choose
the basis so that the essential spectrum of $A$ can be computed {\it in
principle}
from the eigenvalues of the finite dimensional principal submatrices of
$(a_{ij})$?
Roughly speaking, we show that this is possible if the diagonals of the matrix
$(a_{ij})$ tend to zero sufficiently fast.  It is {\it not} possible in
general.
The family of all operators whose matrices satisfy this growth condition is a
Banach $*$-algebra but it is not a \cstar .

In section 4 we seek more precise results which relate to weak$^*$-convergence
as above.  It is appropriate to formulate these results in terms of \cstar s
of operators.  We introduce the notion of {\it degree} of an operator (relative
to a filtration of the underlying Hilbert space into finite dimensional
subspaces), and show that one can expect convergence of this type in the state
space of the given \cstar\ when the \cstar\ contains a dense set of finite
degree
operators and has a unique tracial state.
These results are applied in section 5 to the discretized Hamiltonians
described above.
\subheading{2. Spectral asymptotics of sequences of matrices}

Let $A$ be a bounded self-adjoint operator on a Hilbert space $H$ and suppose
that we have a sequence of finite-dimensional subspaces $H_n$ of $H$ with the
property that the corresponding sequence of projections $P_n\sim H_n$ converges
strongly to the identity.  Normally, the subspaces $H_n$ will be increasing
with $n$, the main examples arising in the cases where we start with an
orthonormal basis $\{e_n: n=1,2,\dots\}$ (resp. $\{e_n: n\in\Bbb Z\}$) and
choose $H_n=[e_1,e_2,\dots,e_n]$ (resp. $H_n=[e_{-n},e_{-n+1},\dots,e_n]\}$).

For every $n\geq 1$, let $A_n$ denote the compression
$$
A_n=P_nA|_{H_n}.\leqno{2.1}
$$
More generally, given an arbitrary sequence $A_1, A_2, \dots$ of finite
dimensional self-adjoint operators,
we introduce two ``spectral" invariants
$\Lambda$, $\Lambda_e$ of the sequence as follows.  $\Lambda$ is defined as the
set of all $\lambda\in\Bbb R$ with the property that there is a sequence
$\lambda_1, \lambda_2,\dots$ such that $\lambda_n\in\sigma(A_n)$ and
$$
\lim_{n\to\infty}\lambda_n=\lambda.
$$
Notice that $\Lambda$ {\it is closed}.  Indeed, a real number $\lambda$ belongs
to the
complement of $\Lambda$ iff
there is an open set $U$ containing $\lambda$ and an infinite sequence of
integers
$n_1 < n_2 < \dots$ with the property that $\sigma(A_{n_k})\cap U = \emptyset$
for
every $k=1,2,\dots$.  Every point of $U$ shares this property with $\lambda$,
and
hence the complement of $\Lambda$ is open.  If the sequence $\{A_n\}$ is
bounded
in norm then $\Lambda$ is compact.
Of course, $\Lambda$ is considerably smaller than the set
$$
L = \cap_n(\overline{\sigma(A_n)\cup\sigma(A_{n+1})\cup\dots})
$$
of limit points of the sequence $\sigma(A_1),\sigma(A_2),\dots$, since limits
along subsequences
may not qualify for membership in $\Lambda$.

For every $n\geq1$ and every
set $S$ of real numbers, let $N_n(S)$ be the number of eigenvalues of $A_n$
which
belong to $S$, counting multiple eigenvalues according to their multiplicity.
$N_n(\cdot)$ is an integer-valued measure which takes values $0, 1,
2,\dots,\dim H_n$.
The points of $\Lambda$ are classified as follows.

\definition{Definition 2.2}
\roster
\item
A point $\lambda\in\Bbb R$ is called essential if, for every open set $U$
containing $\lambda$, we have
$$
\lim_{n\to\infty}N_n(U) = \infty.
$$
The set of essential points is denoted $\Lambda_e$.
\item
$\lambda\in\Bbb R$ is called transient if there is an open set $U$ containing
$\lambda$ such that
$$
\sup_{n\geq1}N_n(U) < \infty.
$$
\endroster
\enddefinition

\remark{Remarks} It is clear that $\Lambda_e\subseteq\Lambda$.
Moreover, notice that $\lambda\in\Bbb R$
is nonessential iff there is an open set $U$ containing $\lambda$ and an
infinite
sequence of integers $n_1<n_2<\dots$ such that
$$
N_{n_k}(U) \leq M < \infty
$$
for every $k=1,2,\dots$.  Thus the nonessential points are an open set in $\Bbb
R$, and
we conclude that {\it the set $\Lambda_e$ of essential points is closed}.
\endremark

It is conceivable that $\Lambda$ may contain anomalous
points which are neither transient nor essential.
In small neighborhoods $U$ of such a point one would find subsequences
$n_k$ for which $N_{n_k}(U)$ remains bounded, and other subsequences $m_k$
for which $N_{m_k}(U)$ is unbounded.
Fortunately, in most situations one has a more manageable dichotomy (see
Theorem 3.8).
At this level of generality, one cannot say much.  We do have the following

\proclaim{Theorem 2.3}Assuming that the sequence $A_1,A_2,\dots$
arises from an operator
$A$ as in 1.1, then we have $\sigma(A)\subseteq\Lambda$ and
$\sigma_e(A)\subseteq\Lambda_e$.
\endproclaim
\demo{proof of $\sigma(A)\subseteq\Lambda$}

Suppose that $\lambda$ is a real
number which does not belong to $\Lambda$.  We will show that $A-\lambda\I$ is
invertible.  Since $\lambda\notin\Lambda$, there is an $\epsilon > 0$
and a subsequence $n_1 < n_2 <  \dots$ such that
$$
\sigma(A_{n_k})\cap(\lambda-\epsilon,\lambda+\epsilon) = \emptyset
$$
for every $k\geq1$.  Thus the distance from $\lambda$ to $\sigma(A_{n_k})$ is
at least
$\epsilon$, hence
$$
\norm(A_{n_k}-\lambda\I)^{-1}\norm \leq 1/\epsilon
$$
for every $k$.  Let
$$
B_k = (A_{n_k} - \lambda\I)^{-1}P_{n_k}.
$$
Since $B_1,B_2,\dots$ is a bounded sequence of self-adjoint operators and
the ball of $\Cal B(H)$ of radius $1/\epsilon$ is weakly sequentially compact,
$\{B_k\}$ has a weakly convergent subsequence.  Thus, by replacing the
original sequence $n_1, n_2,\dots$ with a subsequence of itself, we may assume
that there is a bounded operator $B$ such that
$$
\lim_{k\to\infty}<B_k\xi,\eta> = <B\xi,\eta>,\qquad \xi,\eta\in H.
$$
Since
$$
(A_{n_k}-\lambda\I)B_k = P_{n_k} \leqno{2.4}
$$
for every $k$ and since $(A_{n_k}-\lambda\I)P_{n_k}$ converges to
$A-\lambda\I$ in the {\it strong} operator topology, it follows that the
product $(A_{n_k}-\lambda\I)B_{n_k}$ converges weakly to $(A-\lambda\I)B$.
Hence from 2.4 we may conclude that
$$
(A-\lambda\I)B = \operatornamewithlimits{weak\ lim}_{k\to\infty}P_{n_k} = \I,
$$
proving that $A-\lambda\I$ is right-invertible.  Left-invertiblity follows
by taking adjoints in the preceding equation.
\enddemo

\demo{proof of $\sigma_e(A)\subseteq\Lambda_e$}

Suppose that
$\lambda\notin\Lambda_e$.  We will show that $A-\lambda\I$ is invertible
modulo trace class operators, and hence $\lambda\notin\sigma_e(A)$.
By the hypothesis on $\lambda$, there is a sequence $n_1<n_2<\dots$ of positive
integers and a pair of positive numbers $\epsilon$, $M$ such that
$$
N_{n_k}(\lambda-\epsilon,\lambda+\epsilon) \leq M
$$
for every $k=1,2,\dots$.  For each $k$, consider $A_{n_k}$ to be a self-adjoint
operator in $\Cal B(P_{n_k}H)$ and let $E_k$ be the spectral measure of this
operator.  Then $Q_k = E_k(\lambda-\epsilon,\lambda+\epsilon)$ is a projection
commuting with $A_{n_k}$ satisfying $\roman{dim}Q_k \leq M$, such that the
restriction $B_k$ of $A_{n_k}-\lambda\I$ to the range of $P_{n_k}-Q_k$ is an
invertible operator with $\norm B_k^{-1}\norm \leq 1/\epsilon$.

By using compactness again and passing to a subsequence as in the preceding
argument, we can assume that both sequences $\{B_k^{-1}(P_{n_k}-Q_k)\}$ and
$\{Q_k\}$ converge in the weak operator topology
$$
\align
B_k^{-1}(P_{n_k}-Q_k)&\to C\\
Q_k&\to Q.
\endalign
$$
Notice that $Q$ is a positive trace class operator.  Indeed, since the set
$$
 \{T\in\Cal B(H): T\geq 0, \roman{trace}(T) \leq M \}
$$
is closed in the weak operator topology, the assertion is evident.  Finally,
using the fact that $A_{n_k}-\lambda\I$ converges in the {\it strong} operator
topology to $A-\lambda\I$ and
$$
(A_{n_k}-\lambda\I)B_k^{-1}(P_{n_k}-Q_k) = P_{n_k}-Q_k
$$
for every $k$, we may take weak limits in the latter formula to obtain
$$
(A-\lambda\I)C = \I - Q,
$$
as required\qed
\enddemo

\remark{Remarks} The appendix contains an example which shows that
the inclusions of Theorem 2.3 can be proper.

One can imagine computing the eigenvalues of $A_n$ for every
$n = 1,2,\dots$, picking a point $\lambda\in\Lambda$, and observing the
distribution of eigenvalues in the vicinity of $\lambda$ as $n$ becomes
large.  If $\lambda$ is essential then the number of eigenvalues in any
small interval containing $\lambda$ will increase without limit.  If
$\lambda$ is transient, then there are positive integers $p\leq q$
with the property that eventually you always see at least $p$ points in
the interval but never more than $q$.  Moreover, after $p$ and $q$ are
appropriately adjusted, then the same behavior will occur in {\it every}
sufficiently small interval containing $\lambda$.
\endremark
\subheading{3. Filtrations and their Banach algebras}

Let $A$ be a self adjoint operator on a Hilbert
space $H$.  Suppose we have an orthonormal basis for $H$,
indexed either by $\Bbb N$ or $\Bbb Z$, and we form the
matrix $(a_{ij})$ of $A$ relative to this basis.
If we form a sequence
of finite-dimensional compressions of $(a_{ij})$
then we may compute $\Lambda$ and $\Lambda_e$ in principle, and
the general results of section 2 imply that $\sigma_e(A)\subseteq\Lambda_e$.
In this section we find conditions on the matrix $(a_{ij})$ which
guarantee first, that every point of $\Lambda$ is either essential or
transient and second, that $\sigma_e(A) = \Lambda_e$.  The criterion
is that the series $\sum_k|k|^{1/2}d_k$ should converge, where $d_k$ is the
sup norm of the $k$th diagonal of $(a_{ij})$.  We conclude this
section with some comments about the nature of transient points.

In order to deal with both unilateral and bilateral orthonormal bases
as well as more general situations in which the dimensions of
compressions increase in uneven jumps, it is necessary to work with
filtrations.  We introduce the concept of {\it degree} of an operator
and a Banach $*$-algebra $D(\Cal F)$ of operators associated with a
filtration $\Cal F$.  It is important for the applications that it should be
easy to estimate the norm in $D(\Cal F)$ (see Proposition 3.4).

\definition{Definition 3.1}
\roster
\item
A {\it filtration} of $H$ is a sequence
$\Cal F = \{H_1,H_2,\dots\}$ of finite dimensional subspaces of $H$
such that $H_n\subseteq H_{n+1}$ and
$$
\overline{\cup_n H_n} = H.
$$
\item
Let $\Cal F = \{H_n\}$ be a filtration of $H$ and let $P_n$ be the
projection onto $H_n$.  The {\it degree} of
an operator $A\in\Cal B(H)$ is defined by
$$
\deg (A) = \sup_{n\geq1} \roman{rank}(P_nA-AP_n).
$$
\endroster
\enddefinition

\remark{Remarks}  $\deg (A)$ is either a nonnegative
integer or $+\infty$.  Moreover, we have
\roster
\item"{}"
$\deg(A^*) = \deg(A)$,\quad $\deg(\lambda A) = \deg(A)\ for\ all\
\lambda\neq0$,
\item"{}"
$\deg(A+B)\leq \deg(A)+\deg(B)$,\quad $\deg(AB)\leq \deg(A)+\deg(B)$.
\endroster
Only the last of the four properties is not quite obvious; it follows
from the fact that the map $A\to [P_n,A] = P_nA-AP_n$ is a derivation,
$[P_n,AB] = [P_n,A]B + A[P_n,B]$, which implies
$$
\roman{rank}[P_n,AB]\leq \roman{rank}[P_n,A] + \roman{rank}[P_n,B].
$$
We conclude from these observations that the set of all finite-degree
operators in $\Cal B(H)$
is a self-adjoint unital subalgebra of $\Cal B(H)$.
\endremark

\remark{Remark}  Notice that for any operator $A$ and any projection
$P$ we have
$$
PA-AP = PA(\I - P) - (\I - P)AP,
$$
hence
$$
\align
\roman{rank}(PA - AP) &= \roman{rank}(PA(\I - P)) + \roman{rank}((\I - P)AP) \\
                      &= \roman{rank}((\I - P)A^*P) + \roman{rank}((\I - P)AP).
\endalign
$$
It follows that $\roman{deg}(A) < \infty$ iff both operators $B = A$ and
$B=A^*$ satisfy the condition
$$
\sup_{n\geq1} \roman{rank}((\I - P_n)BP_n) < +\infty.
$$
Thus, {\it operators of finite degree are abstractions of band-limited
matrices}.
\endremark

We associate a Banach $*$-algebra $D(\Cal F)$ to a filtration $\Cal F$ in the
following
way.  Let $D(\Cal F)$ denote the set of all operators $A\in\Cal B(H)$
having a decomposition into an infinite sum of finite degree operators $A_k$
$$
A = \sum_{k=1}^\infty A_k \tag{3.2}
$$
in such a way that the sum
$$
s=\sum_{k=1}^\infty (1+\roman{deg}(A_k)^{1/2})\|A_k\|
$$
is finite.  Notice that finiteness of $s$ ensures that the series
3.2 converges absolutely with respect to the operator norm.  We define
$|A|_\Cal F$ to be the infimum of all such sums $s$ which arise from
representations
of $A$ as in 3.2.

\proclaim{Proposition 3.3}  With respect to the norm $|\cdot|_\Cal F$ and the
operator adjoint, $D(\Cal F)$ is a unital Banach $*$-algebra
in which $|\I|_\Cal F =1$.
\endproclaim
\demo{proof}  Let $\Cal D$ be the $*$-algebra of all finite degree operators in
$\Cal B(H)$, and consider the function $\phi:\Cal D\to\Bbb R^+$ defined by
$$
\phi(A) = (1+(\roman{deg}A)^{1/2})\|A\|.
$$
The norm on $D(\Cal F)$ is defined by
$$
|A|_\Cal F = \inf \sum_{k=1}^{\infty}\phi(A_k),
$$
the infimum extended over all sequences $A_1,A_2,\dots$ in $\Cal D$ satisfying
$$
A = \sum_{k=1}^\infty A_k
$$
together with $\sum_k \phi(A_k) < +\infty$.  The fact that $(D(\Cal
F),|\cdot|_\Cal F)$ is
a Banach $*$-algebra follows from straightforward applications of the following
properties of $\phi$:
\roster
\item
$\phi(A) \geq\|A\|.$
\item
$\phi(\lambda A) = |\lambda|\cdot \|A\|, \quad \lambda\in\Bbb C.$
\item
$\phi(A^*) = \phi(A).$
\item
$\phi(AB) \leq\phi(A)\phi(B).$
\endroster
For example, (4) follows from the fact that the weight sequence
$w_0,w_1,w_2,\dots$ defined by $w_k=1 + k^{1/2}$ satisfies
$$
w_{k+j} \leq w_k w_j, \qquad k,j\geq 0,
$$
together with
$$
\roman{deg}(AB) \leq \roman{deg}A + \roman{deg}B.
$$
The identity of $\Cal B(H)$ is of degree 0, hence $|\I|_\Cal F \leq 1$, while
the
inequality $|\I|_\Cal F\geq1$ is valid for the unit of any Banach algebra
\qed\enddemo

$D(\Cal F)$ is certainly not a \cstar, but it is dense in $\Cal B(H)$ in the
strong operator topology.  The following gives a concrete description of
operators in $D(\Cal F)$ for one of the two most important filtrations.

\proclaim{Proposition 3.4}   Let $\{e_n: n\in \Bbb Z\}$ be a bilateral
orthonormal
basis for a Hilbert space $H$ and let $\Cal F=\{H_1, H_2,\dots\}$
be the filtration defined by
$$
H_n = [e_{-n},e_{-n+1},\dots,e_n].
$$
Let $(a_{ij})$ be the matrix of an operator
$A\in\Cal B(H)$ relative to $\{e_n\}$, and for every $k\in\Bbb Z$ let
$$
d_k=\sup_{i\in\Bbb Z} |a_{i+k,i}|
$$
be the sup norm of the $k$th diagonal of $(a_{ij})$.  Then
$$
|A|_\Cal F \leq \sum_{k=-\infty}^{+\infty}(1 + |2k|^{1/2})d_k.\tag{3.5}
$$
In particular, $A$ will belong to $D(\Cal F)$ whenever the series
$\sum_k |k|^{1/2}d_k$ converges.
\endproclaim

\demo{proof}
We may assume that the sum on the right side of 3.5 converges.  Let $D_k$
be the operator whose
matrix agrees with $(a_{ij})$ along the kth diagonal and is zero elsewhere.
Notice that $||D_k|| = d_k$, and hence
$$
A = \sum_{k=-\infty}^{+\infty} D_k,
$$
the sum on the right converging absolutely in the operator norm because
$$
\sum_k |k|^{1/2}d_k <\infty.
$$

We claim that $\roman{deg}D_k\leq 2|k|$ for every $k\in\Bbb Z$.  To see that,
fix $n=1,2,\dots$ and let $P_n$ be the projection onto $H_n$.
The operator $D_k$ maps the range of
$P_n$ into the range of $P_{n+|k|}$, and a careful inspection of the two cases
$k>0$ and $k<0$ shows that the dimension of $(\I-P_n)P_{n+|k|}$ is
at most $|k|$ for every $n\geq1$.  It follows that the rank of either operator
$(\I-P_n)D_kP_n$ or $(\I-P_n)D_k^*P_n$ is at most $|k|$, independently of $n$.
Thus from a previous remark we conclude that $\roman{deg}D_k\leq 2|k|$.

Now $\|D_k\| = d_k$, and since the series   $\sum_k |k|^{1/2}d_k$ converges
we must have $\sum_k \|D_k\| < \infty$.  Hence
$$
A = \sum_{k=-\infty}^{+\infty}D_k
$$
is an absolutely convergent sum of operators $D_k$ satisfying
$\roman{deg}D_k \leq 2|k|$, from which 3.5 follows\qed
\enddemo

In particular, any operator whose matrix $(a_{ij})$ is {\it band-limited} in
the sense that $a_{ij} = 0$ whenever $|i-j|$ is sufficiently large, must
belong to $D(\Cal F)$.

The following estimate occupies a key position, and it explains why the
norm on $D(\Cal F)$ is defined as it is.  For $B\in\Cal B(H)$, we write
$\|B\|_2 = (\trace(B^*B))^{1/2}$ for the Hilbert-Schmidt norm of B.

\proclaim{Lemma 3.6} Let $\{H_1, H_2,\dots\}$ be a filtration of $H$ and
let $P_n$ be the projection onto $H_n$.  Then for every $A\in D(\Cal F)$ we
have
$$
\sup_{n\geq1}\|AP_n - P_nAP_n\|_2 \leq |A|_\Cal F < \infty.
$$
\endproclaim
\demo{proof}
We claim first that for every $A\in\Cal B(H)$ and every $n=1,2,\dots$ we have
$$
\trace|AP_n - P_nAP_n|^2 \leq \roman{deg}(A) \|A\|^2,\tag{3.7}
$$
where for an operator $B\in \Cal B(H)$, $|B|$ denotes the positive square
root of $B^*B$.  Indeed, since $\trace |B|^2\leq \roman{rank}(B) \|B\|^2$,
we can estimate the left side of 3.7 as follows
$$
\align
\trace|(\I - P_n)AP_n|^2 &\leq \roman{rank}((\I-P_n)AP_n)) \|(\I-P_n)AP_n\|^2\\
   &= \roman{rank}((\I-P_n)(AP_n-P_nA)) \|(\I-P_n)AP_n\|^2 \\
   &\leq \roman{rank}(AP_n-P_nA) \|A\|^2 \leq \roman{deg}(A) \|A\|^2.
\endalign
$$
By taking the square root of both sides of 3.7 we obtain
$$
\|AP_n-P_nAP_n\|_2 \leq \roman{deg}(A)^{1/2} \|A\|.
$$

Now suppose $A\in D(\Cal F)$, and choose finite degree operators $A_k$ such
that
$$\sum_k(1 + \roman{deg}(A_k)^{1/2}) \|A_k\| < \infty$$
and $A = A_1 + A_2 + \dots$.
Using the triangle inequality for the Hilbert-Schmidt norm we see from 3.7 that
$$
\align
\|P_nA - P_nAP_n\|_2 &\leq \sum_{k=1}^\infty\roman{deg}(A_k)^{1/2}\|A_k\| \\
 &\leq \sum_{k=1}^\infty (1+\roman{deg}(A_k)^{1/2}) \|A_k\|.
\endalign
$$
The desired inequality follows by taking the supremum over $n$ and the infimum
over all such sequences $A_1, A_2,\dots$ \qed
\enddemo

\remark{Remark}  Notice that Lemma 3.6 implies that for every $A\in D(\Cal F)$
we have the following uniform estimate of commutator $2$-norms
$$
\sup_n\|AP_n - P_nA\|_2 \leq \sqrt2 |A|_\Cal F < +\infty.
$$
\endremark

Now let $A$ be a self-adjoint operator in $\Cal B(H)$, let
$\Cal F = \{H_1,H_2,\dots\}$ be a filtration, and let $A_n$ be the
compression of $A$ to $H_n$.
The following result is our basic description of the
spectrum of $A$ in terms of the eigenvalues of $\{A_n\}$.
It asserts that the essential spectrum of $A$ is precisely
the set of essential points, and that all of the remaining points
of $\Lambda$ are transient.

\proclaim{Theorem 3.8} Assume that $A = A^*$ belongs to the Banach
algebra $D(\Cal F)$.  Then
\roster
\item"{(i)}"
$\sigma_e(A) = \Lambda_e$.
\item"{(ii)}"
Every point of $\Lambda$ is either transient or essential.
\endroster
\endproclaim

\demo{proof} We will prove that every point in the complement of
the essential spectrum of $A$ is transient.  Notice
that both statements (i) and (ii) follow from this.  Indeed, because
of Theorem 2.3 we know that $\sigma_e(A)\subseteq\Lambda_e$.  On the
other hand, since no transient point can be an essential point, the
above assertion implies that the complement of $\sigma_e(A)$ is contained
in the complement of $\Lambda_e$; hence (i).  Armed with (i), we see that
$\Lambda\setminus\Lambda_e$ is the same as $\Lambda\setminus\sigma_e(A)$,
which by the assertion is contained in the transient points of $\Lambda$.
Thus we obtain (ii).

To prove this assertion, choose a point
$\lambda$ in the complement of $\sigma_e(A)$.  In order to show that
$\lambda$ is transient,
we claim first that there is a finite-dimensional projection $Q$ and a
positive number $\epsilon$ such that for every real number
$\mu\in[\lambda-\epsilon,\lambda+\epsilon]$, the operator
$$
A+Q-\mu\I
$$
is invertible.  Indeed, since $\lambda\notin\sigma_e(A)$ it follows that
either $\lambda$ is an isolated point of the spectrum of $A$ whose
eigenspace is finite dimensional, or else $A-\lambda\I$ is invertible.
In either case
$$
\{\xi\in H: A\xi=\lambda\xi\}
$$
is finite dimensional and if $Q$ denotes the projection onto
this eigenspace, then $Q$ commutes with $A$ and $A+Q-\lambda\I$
is an invertible self-adjoint
operator.  Because the invertible operators are an open set
it follows that there is a
neighborhood $U$ of $\lambda$ with the property that $A+Q-\mu\I$ is
invertible for every $\mu\in U$.  We can take
$[\lambda-\epsilon,\lambda+\epsilon]$ to be an appropriate closed interval
about $\lambda$ which is contained in such a neighborhood.

Note too that by continuity of inversion, there is a positive number $M$
such that
$$
\norm (A+Q-\mu\I)^{-1}\norm  \leq M\tag{3.9}
$$
for every $\mu$ satisfying $|\mu-\lambda|\leq\epsilon$.

We will
show that if $V$ is the open interval $(\lambda-\epsilon,\lambda+\epsilon)$
and $N_n(V)$ is the number of eigenvalues of $A_n$ which lie in $V$, then
the sequence $N_n(V)$ is bounded.  Thus
$\lambda$ must be a transient point, and that will complete the proof.
To that end, note first that by Lemma 3.6 there is a positive number $N$
(in fact, $N = |A|_{\Cal F}^2$ will do) such that
$$
\trace |AP_n-P_nAP_n|^2\leq N,
$$
for every $n=1,2,\dots$.  Fix $n$, and let
$\lambda_1,\lambda_2,\dots,\lambda_p$
be the set of all eigenvalues of $A_n$ which lie in $V$, with repetitions
according to the multiplicity of the corresponding eigenspace.  Actually, we
will show that
$$
p\leq M^2N + \dim Q. \tag{3.10}
$$
Since the right side of 3.10 does not depend on $n$ it will follow that
$\lambda$ is transient, completing the proof.

To prove 3.10, we choose for each $k=1,2,\dots,p$ a unit vector
$\xi_k$ in the domain of $A_n$ (i.e., the range of $P_n$) such that
$A_n\xi_k = \lambda_k\xi_k$, and such that $\{\xi_1,\dots,\xi_p\}$
is an orthonormal set.  This is possible because the eigenvectors
which belong to different eigenvalues of a self-adjoint operator are mutually
orthogonal; in the case of a multiple eigenvalue we simply choose an
orthonormal basis for that eigenspace.

Note first that for each $k$ between $1$ and $p$,
$$
A\xi_k-\lambda_k\xi_k = (\I - P_n)A\xi_k.    \tag{3.11}
$$
Indeed, since $P_n(A\xi_k - \lambda_k\xi_k)= A_n\xi_k - \lambda_k\xi_k=0$,
we can write
$$
A\xi_k - \lambda_k\xi_k = (\I - P_n)(A\xi_k - \lambda_k\xi_k)
     =(\I - P_n)A\xi_k.
$$
We claim now that for each $k$ between $1$ and $p$, we have the inequality
$$
1 \leq M^2\norm (\I - P_n)A\xi_k\norm ^2 + <Q\xi_k,\xi_k>.
 \tag{3.12}
$$
To see that, write
$$
\align
1 - <Q\xi_k,\xi_k>
&=\ \norm(\I - Q)\xi_k\norm^2 \ \leq\  M^2
\norm(A+Q-\lambda_k)(\I-Q)\xi_k\norm^2 \\
&= M^2\norm(A-\lambda_k\I)(\I -Q)\xi_k\norm^2\ =\ \norm(\I -Q)(A\xi_k -
\lambda_k\xi_k)\norm \\
&= M^2 \norm(\I - Q)(\I - P_n)A\xi_k\norm^2\ \leq\  M^2 \norm(\I
-P_n)A\xi_k\norm^2.
\endalign
$$
The first inequality follows from 3.9, the identity
$(A - \lambda_k)(\I -Q)\xi_k = (\I -Q)(A\xi_k - \lambda_k\xi_k)$
follows from the fact that $Q$ and $A$ commute, and the last equality
and inequality follow respectively from 3.11 and
the fact that $\I - Q$ is a contraction.
Summing the inequality 3.12 on $k$ we obtain
$$
\align
p &\leq M^2\sum_{k=1}^p\norm(\I - P_n)A\xi_k\norm^2 + \sum_{k=1}^p
<Q\xi_k,\xi_k> \\
  &\leq M^2\,\trace|(\I -P_n)AP_n|^2 + \trace(Q) \\
  &\leq M^2 N + \roman{dim}Q,
\endalign
$$
and 3.10 is established \qed
\enddemo

\remark{Remark 3.13} We conclude this section with some
remarks about transient points and other spurious eigenvalues.
To provide a context
for these remarks, let us return to
the classical setting of Szeg\"o's theorem, in which $A_1,A_2,\dots$ is a
sequence of Toeplitz matrices.  Let $f$ be a bounded
measurable real-valued function on the interval $[-\pi,\pi]$, and consider
the Fourier series of $f$
$$
f(x) \sim \sum_{n=-\infty}^\infty a_n e^{inx}.
$$
For $n = 1,2,\dots$ let $A_n$ be the Toeplitz matrix
$$
A_n = \pmatrix
a_0&a_{-1}&\hdots&a_{-n+1}\\
a_1&a_0&\hdots&a_{-n+2}\\
\vdots&&\ddots&\vdots\\
a_{n-1}&\hdots&&a_0
\endpmatrix.
$$
Let $a$ and $b$ be respectively the essential inf and sup of $f$
over $[\pi,\pi]$.  It is known that the eigenvalues of the sequence
$\{A_n: n=1,2,\dots\}$ tend to fill out the entire interval $[a,b]$
(see \cite{11, pp 201--202} and \cite{12, Lemma II.1}).
In more concrete terms, suppose that $\lambda\in (a,b)$ and that we are looking
in a small neighborhood $U=(\lambda - \epsilon,\lambda + \epsilon)$ of
$\lambda$.  Then $U\cap\sigma(A_n)\neq\emptyset$ for sufficiently large $n$,
regardless of where $\lambda$ is located in the interval $(a,b)$.  The
point we want to make is that this can be misleading, and that what one is
observing here may have nothing to do with the spectrum of the operator $A$.

Indeed, for these operators $A$ we have
$$
\sigma(A)=\sigma_e(A)=R,
$$
where $R$ is the essential range of the function $f$.  Suppose that $f$
is chosen so that its essential range is disconnected.  If we choose
a point $\lambda\in[a,b]\setminus R$
then for every $\epsilon>0$ we will find the sets of eigenvalues
$(\lambda-\epsilon,\lambda+\epsilon)\cap\sigma(A_n)$ to be permanently nonempty
for large $n$, in spite of the fact that $(\lambda-\epsilon,\lambda+\epsilon)$
contains no points of the spectrum of $A$ if $\epsilon$ is small.

Returning to the context of Theorem 3.8, it is not hard to give examples of
band-limited matrices $(a_{ij})$ which exhibit the same kind of behavior,
namely
$\sigma(A) = \sigma_e(A)$ while $\Lambda$ contains points not in
the spectrum of the operator.  In this case, 3.8 implies that such points will
be
transient.  More generally, if one starts with a band-limited
operator whose spectrum is
properly larger than its essential spectrum, then every point of
$\sigma(A)\setminus\sigma_e(A)$ will be a transient point (Theorem 3.8).
Of course, for such an operator there may also be other transient points which
do {\it not} belong to $\sigma(A)$.  Unfortunately, we do not know how to
distinguish between transient points which belong to the spectrum and
transient points which do not.  Perhaps it is not even possible to do so.
In any case, these observations have led us to the conclusion that
{\it transient points should be ignored}.  In doing that one is also ignoring
points of $\sigma(A)\setminus\sigma_e(A)$; but since such points are merely
isolated eigenvalues having finite multiplicity, the cost is small.

We conclude that if one is interested in approximating
spectra of infinite dimensional operators using matrix truncations,
then one should not simply look for the eigenvalues of the truncations.
Rather, {\it one should weight the
count of eigenvalues in a way which eliminates spurious ones occuring
in the vicinity of transient points}.
The simplest example of such weighting is the sequence of `densities'
$$
{N_n(U)\over n},
$$
where $N_n(U)$ is the number of eigenvalues of the $n\times n$ truncation
of $A$ which belong to the open set $U$.
Indeed, in the present context we obviously have
$$
\lim_{n\to\infty}{N_n(U)\over n} = 0
$$
whenever $U$ is a sufficiently small neighborhood of a transient point
$\lambda$, simply because the sequence $N_1(U), N_2(U),\dots$ is bounded.
If $\lambda$ is an essential point on the other hand, then we require
conditions under which the limit
$$
\lim_{n\to\infty}{N_n(U)\over n}
$$
exists and is positive for every open neighborhood $U$ of $\lambda$.
The existence of this limit in the case of the multiplication
operators $A$ discussed above is a consequence of
Szeg\"o's theorem, as we have pointed out in the introduction.
In the following section we take up the asymptotic behavior of
these densities in a more general context.
\endremark

\subheading{4. Filtrations and convergence of densities}

In this section we shift our point of view somewhat. Rather than
consider single self-adjoint operators we consider a concretely
presented \cstar\  of operators $\Cal A \subseteq\Cal B(H)$ and we
introduce the concept of an $\Cal A$-filtration.  This
is an abstraction of an orthonormal
basis for $H$ with respect to which a dense set of the operators
in $\Cal A$ are band-limited.  The point we want to make
is that in order to use the methods of this paper to
compute the spectrum of a self-adjoint operator in a
\cstar\ $\Cal C$, one
should first look for a faithful representation
$\pi: \Cal C\to \Cal B(H)$ with the property that the \cstar\
$\Cal  A=\pi(\Cal C)$
admits a natural $\Cal A$-filtration, and then proceed with the analysis
associated with Theorem 3.8 or Theorem 4.5 below.

After presenting some examples, we prove
a counterpart of Szeg\"o's theorem which is
appropriate for \cstar s.  In order
to keep the statement and proof as simple as possible, we
assume that the \cstar\ has a unique trace.  While this is not
the most general formulation possible, it does provide the
basis for the applications which will be taken up in the next
section.  It will be convenient (though not necessary) to assume
that \cstar s of operators on $H$ contain the identity of $\Cal B(H)$.

\definition{Definition 4.1}  Let $\Cal A\subseteq\Cal B(H)$ be a
\cstar.  An $\Cal A$-filtration is a filtration of $H$
with the property that the $*$-subalgebra of all finite degree
operators in $\Cal A$ is norm-dense in $\Cal A$.
\enddefinition

\example{Example 1}  The simplest example is that in which $\Cal A$ is an
AF-algebra
having a cyclic vector $\xi$.  If $\Cal A_1\subseteq\Cal A_2\subseteq\dots$
is a sequence of finite dimensional *-subalgebras of
$\Cal A$ whose union is norm-dense in $\Cal A$, then $H_n = [\Cal A_n\xi]$
defines a
filtration for which every element of the union $\cup_n\Cal A_n$ has finite
degree.
\endexample

\example{Example 2}
Suppose on the other hand that we have an action of
$\Bbb Z$ on a compact Hausdorff space $X$
for which there is a point $x_0\in X$ whose orbit under
$\Bbb Z$ is dense.  Then one can write down a
representation $\pi$ of the \cstar ic crossed product
$$
\Cal C=\Bbb Z\times C(X)
$$
and a natural filtration for the \cstar\ $\Cal A=\pi(\Cal C)$.  To see that,
let $\{e_n: n\in\Bbb Z\}$ be an orthonormal basis for a Hilbert space $H$.
The representation $\pi_0:C(X)\to \Cal B(H)$ and the unitary operator
$U\in\Cal B(H)$ defined by
$$
\align
\pi_0(f)e_n&=f(n\cdot x_0)e_n\\
Ue_n&=e_{n+1}
\endalign
$$
define a covariant pair for the induced action of $\Bbb Z$ on $C(X)$, and
hence the pair $(\pi_0,U)$ defines a representation
$\pi$ of the crossed product.  Let
$$
H_n=[e_{-n},\dots,e_n], \qquad n\geq1,
$$
and let $\Cal D$ be the image of
$C(X)$ under $\pi_0$.  It is clear that every element of $\Cal D$
has degree zero, and that $\roman{deg}(U) = 1$.
Thus any operator expressible as a finite sum
$$
A = \sum_{k=-n}^{n}D_kU^k
$$
where $D_k\in\Cal D$, has finite degree.  Since these
operators are dense in $\Cal A$ it follows that
$\{H_1\subseteq H_2\subseteq\dots\}$ is an $\Cal A$-filtration.
\endexample

\example{Example 3}
Let $\Cal T$ be the Toeplitz \cstar.  We represent $\Cal T$ as the
\cstar\ generated by a simple unilateral shift $S$, acting on an
orthonormal basis $\{e_1,e_2,\dots\}$ for $H$ via
$Se_n=e_{n+1}$.  Let $\Cal F=\{H_n: n\geq1\}$ be the filtration
$H_n=[e_1,e_2,\dots,e_n]$.  Notice that $\roman{deg}(S)=1$.  Indeed,
letting $P_n$ be the projection onto $H_n$, we have
$$
SP_n=P_nSP_n + (P_{n+1} - P_n)SP_n,
$$
while $P_nS=P_nSP_n$.  Hence $SP_n - P_nS = (P_{n+1} - P_n)SP_n$
is a rank-one operator for every $n$.  It follows that the $*$-algebra
$\Cal T_0$ generated by $S$ consists of finite degree operators, and
$\Cal T_0$ is dense in $\Cal T$.  Thus $\Cal F$ is a $\Cal T$-filtration.
\endexample

We require the following estimate.

\proclaim{Proposition 4.2}  Let $\Cal F=\{H_n\}$ be a filtration of $H$,
let $P_n$ be the projection onto $H_n$, and let $A_1,A_2,\dots,A_p$ be a
finite set of operators in $\Cal B(H)$.
Then for every $n=1,2,\dots$ we have
$$
\trace|P_nA_1A_2\dots A_pP_n - P_nA_1P_nA_2P_n\dots P_nA_pP_n|
\leq \|A_1\|\dots\|A_p\| \sum_{k=1}^p\deg A_k.
$$
\endproclaim

\demo{proof}
Without loss of generality, we can assume that each operator $A_k$ is a
contraction and has finite degree.  Let us fix $n$, and define
$$
\Delta_p = P_nA_1A_2\dots A_pP_n - P_nA_1P_nA_2P_n\dots P_nA_pP_n.
$$
We claim first that there are contractions
$B_1,B_2,\dots,B_p, C_1, C_2,\dots,C_p$ such that
$$
\Delta_p = \sum_{k=1}^p B_k(\I-P_n)A_kP_nC_k.
$$
This is easily seen by induction on $p$.  It is trivial for
$p=1$ (take $B=C=0$).  Assuming it valid for $p$, we have
$$
\Delta_{p+1} = P_1A_1\dots A_p(\I-P_n)A_{p+1}P_n + \Delta_pA_{p+1}P_n,
$$
and the conclusion is evident from the induction hypothesis on $\Delta_p$.

Each operator $(\I - P_n)A_kP_n$ is of finite rank, since
$$
\align
\roman{rank}((\I - P_n)A_pP_n) &=\roman{rank}((\I - P_n)(A_pP_n - P_nA_p))\\
&\leq \roman{rank}(A_pP_n - P_nA_p) \leq \roman{deg} A_p.
\endalign
$$
Thus we may write

$$
\align
\trace|\Delta_p|&\leq \sum_{k=1}^p \trace|B_k(\I-P_n)A_kP_nC_k| \\
&\leq \sum_{k=1}^p\|B_k\| \|C_k\|\, \trace|(\I-P_n)A_kP_n| \tag{4.3}\\
&\leq \sum_{k=1}^p\trace|(\I-P_n)A_kP_n| ,
\endalign
$$
where we have used the elementary fact that for any finite rank
operator $F$,
$$
\trace|BFC|\leq \|B\| \|C\|\,\trace|F|.
$$
Moreover, for such an $F$ we have
$$
\trace|F| \leq \|F\|\cdot \roman{rank}(F).
$$
Hence
$$
\align
\trace|(\I-P_n)A_kP_n| &\leq \roman{rank}((\I-P_n)A_kP_n) =
\roman{rank}((\I-P_n)(A_kP_n-P_nA_k))\\
                       &\leq \roman{rank}(A_kP_n-P_nA_k) \leq \deg A_k.
\endalign
$$
Utilizing the latter inequality in (4.3) we obtain
$$
\trace|\Delta_p|\leq \sum_{k=1}^p\deg A_k,
$$
as required.
\qed
\enddemo

Let $\Cal A$ be a concrete \cstar.  A tracial state of $\Cal A$ is a
positive linear functional $\rho$ satisfying $\rho(\I)=1$ and
$\rho(AB)=\rho(BA)$ for every $A,B\in \Cal A$.  One might try
to construct tracial states of $\Cal A$ by using a filtration
$\Cal F = \{H_n: n\geq1\}$ for $H$ in the following way.  Letting
$P_n$ be the projection on $H_n$, we can define a state $\rho_n$ of
$\Cal B(H)$ by way of
$$
\rho_n(T) = {1\over {\roman{dim}(H_n)}}\trace(P_nT).
$$
The restriction of $\rho_n$ to $P_n\Cal B(H)P_n$ is a trace.  Of course,
$\rho_n$ does not restrict to a trace of $\Cal A$, but one might attempt
to obtain a trace by taking weak*-limits of subsequences of
$\{\rho_n\restriction_{\Cal A}: n\geq1\}$, or their averages.  In general,
this program will fail.  It will certainly fail if $\Cal A$ is
isomorphic to a Cuntz algebra $\Cal O_n, n=2,3,\dots,\infty$, since
$\Cal O_n$ has no tracial states whatsoever.  On the other hand,
the following result
shows that this procedure will succeed when $\Cal F$ is an
$\Cal A$-filtration.

\proclaim{Proposition 4.4}  Suppose that $\{H_1\subseteq H_2\subseteq\dots\}$
is an $\Cal A$-filtration.  For every $n=1,2,\dots$ put $d_n=\roman{dim}H_n$,
and
let $\rho_n$ be the state of $\Cal A$ defined by
$$
\rho_n(A) = {1\over {d_n}}\trace(P_nA).
$$
Let $R_n$ be the weak$^*$-closed convex hull of the set
$\{\rho_n,\rho_{n+1},\rho_{n+2},\dots\}$.  Then $\cap_n R_n$ is a
nonempty set of tracial states.
\endproclaim
\demo{proof}
$R_1,R_2,\dots$ is a decreasing sequence of nonvoid compact convex subsets
of the state space of $\Cal A$, and hence the intersection $R_\infty$ is a
nonvoid compact convex set of states.  We have to show that every element
$\sigma\in R_\infty$ is a trace.  Since the finite degree elements are
dense in $\Cal A$ and $\sigma$ is bounded, it suffices to show that
$$
\sigma(AB) = \sigma(BA)
$$
for operators $A,B\in\Cal A$ having finite degree.  Because of the
inequality
$$
|\sigma(AB)-\sigma(BA)|\leq \limsup_{n\to\infty}|\rho_n(AB) - \rho_n(BA)|,
$$
it suffices to show that $|\rho_n(AB) - \rho_n(BA)|\to 0$ as $n\to\infty$.
But since
$$
\trace(P_nAP_nBP_n) = \trace(P_nBP_nAP_n)
$$
for every $n$, we may
use the case $p=2$ of Proposition 4.2 as follows
$$
\align
|\rho_n(AB) &- \rho_n(BA)|=
{1\over{d_n}}|(\trace(P_nABP_n)-\trace(P_nBAP_n)|\\
&\leq{1\over{d_n}}(|\trace(P_nABP_n-P_nAP_nBP_n)| +
|\trace(P_nBAP_n-P_nBP_nAP_n)|)\\
&\leq {2\over{d_n}}\|A\|\|B\|(\roman{deg}A+ \roman{deg}B).
\endalign
$$
The right side obviously tends to $0$ as $n\to\infty$\qed
\enddemo

\remark{Remarks}  Given a \cstar\ $\Cal A\subseteq\Cal B(H)$, it is
reasonable to ask if it is always possible to find an $\Cal A$-filtration.
Proposition 4.4 implies that {\it if}
an $\Cal A$-filtration exists then $\Cal A$ must have at least
one tracial state.  Thus, many concrete \cstar s cannot be
associated with a
filtration.  For example, if $\Cal A$ is isomorphic to one of
the Cuntz algebras $\Cal O_n$, $n=2,3,\dots,\infty$, then
$\Cal A$ has no tracial states whatsoever, and hence
$\Cal A$-filtrations do not exist.
\endremark

We assume throughout the remainder of this section
that $\Cal A$ has a {\it unique} tracial state $\tau$.
Then every self-adjoint operator $A\in\Cal A$ determines
a natural probability measure $\mu_A$ on $\Bbb R$ by way of
$$
\int_{-\infty}^{+\infty}f(x)\,d\mu_A(x) = \tau(f(A))
$$
for every $f\in C_0(\Bbb R)$.  Clearly the closed support of $\mu_A$
is contained in the
spectrum of $A$, and if $\tau$ is a faithful trace then we have
$$
\roman{support}(\mu_A) = \sigma(A).
$$
Here is a more explicit description of $\mu_A$.
The GNS construction applied to $\tau$ gives rise to
a representation of $\Cal A$ on another Hilbert space in which
the image of $\Cal A$ generates a finite factor $M$.
The spectral measure $E_A$ of $A$ takes values in the projections of $M$.
Letting $\roman{tr}$ denote the natural normalized trace on $M$ we have
$$
\mu_A(S) = \roman{tr}(E_A(S))
$$
for every Borel set $S\subseteq\Bbb R$.  The measure $\mu_A$ will be called
the {\it spectral distribution} of $A$.

Assuming now that $\Cal F=\{H_n\}$ is an $\Cal A$-filtration, we have

\proclaim{Theorem 4.5} Let $\mu_A$ be the
spectral distribution of a self-adjoint operator $A\in\Cal A$, and
let $[a,b]$ be the smallest closed interval containing $\sigma(A)$.
For each $n$, let $d_n=\dim H_n$ and
let $\lambda_1,\lambda_2,\dots,\lambda_{d_n}$ be the eigenvalues of
$A_n=P_nA\restriction H_n$, repeated according to multiplicity.  Then
for every $f\in C[a,b]$,
$$
\lim_{n\to\infty}{1\over
{d_n}}(f(\lambda_1)+f(\lambda_2)+\dots+f(\lambda_{d_n}))
= \int_a^b f(x)\, d\mu_A(x).
$$
\endproclaim
\demo{proof}
Let $\tau_n$ be the state of $\Cal B(H)$ defined by
$$
\tau_n(T) = {1\over d_n}\trace (P_nT).
$$
Noting that $\tau_n$ restricts to the normalized trace on $P_n\Cal B(H)P_n$,
Theorem 4.5 will become evident once we have established the following two
assertions.
\roster
\item"{(i)}" For every $f\in C[a,b]$, we have
$$
\lim_{n\to\infty}|\tau_n(f(A))-\tau_n(f(P_nAP_n))| = 0.
$$
\item"{(ii)}" For every $B\in\Cal A$,
$$
\lim_{n\to\infty}\tau_n(B) = \tau(B),
$$
where $\tau$ is the tracial state of $\Cal A$.
\endroster

To prove (i), notice that since the sequence of linear functionals
$$
f\in C[a,b]\mapsto \tau_n(f(A)) - \tau_n(f(P_nAP_n))
$$
is uniformly bounded (an upper bound for their norms is 2), it is
enough to prove (i) for polynomials $f$; and by linearity, we may
further reduce to the case where $f$ is a monomial, $f(x)=x^p$ for
some $p=1,2,\dots$.  Finally, since the sequence of $p$-linear
forms
$$
B_n(T_1,T_2,\dots,T_p) =
     \tau_n(T_1T_2\dots T_p)-\tau_n(P_nT_1P_nT_2P_n\dots P_nT_pP_n)
$$
is uniformly bounded (again, $\|B_n\|\leq 2$), it is enough to prove that
$$
\lim_{n\to\infty}|\tau_n(A^p) - \tau_n((P_nAP_n)^p)|=0
$$
for operators $A$ having finite degree.  But in this case, Proposition 4.2
implies
that
$$
|\tau_n(A^p) - \tau_n((P_nAP_n)^p)|
\leq {p\over {d_n}} \|A\|^p \roman{deg}(A),
$$
and the right side obviously tends to 0 as $n\to\infty$.

Assertion (ii) follows immediately from Proposition 4.4.  For if $R_n$ denotes
the sequence of convex sets of 4.5, then their intersection must consist of
the singleton $\{\tau\}$; hence the sequence of states
$$
B\in\Cal A\mapsto\tau_n(B)= {1\over{d_n}}\trace(P_nB)
$$
has a unique weak$^*$ cluster point $\tau$, hence it must actually
converge to $\tau$ in the weak$^*$ topology
\qed
\enddemo

\remark{Remark 4.6}  Let $A=A^*\in \Cal A$.
If one drops the hypothesis that $\Cal A$ has a unique
trace, then the averages of the eigenvalues
of $P_nA\restriction{H_n}$ may fail to converge.  However, one can push the
method
of proof to obtain other useful conclusions in certain cases.  For example,
if we assume that there are {\it sufficiently many} traces in the sense
that for every nonzero $B\in\Cal A$ there is a tracial state $\tau$
such that $\tau(B^*B)>0$, then one can obtain the following information
about the density of eigenvalues in the vicinity of a point of the spectrum
of $A$.  Let $N_n(U)$ denote the number of eigenvalues of
$P_nA\restriction{H_n}$
which belong to the open set $U$, and let $\lambda$ be a point in the spectrum
of $A$.  Then for every neighborhood $U$ of $\lambda$ there is a positive
number $\epsilon$ such that
$$
\frac{N_n(U)}{\roman{dim}(H_n)}\geq \epsilon,
$$
for all sufficiently large $n$.  We omit the details since this result is
not relevant to our needs below.
\endremark
\subheading{5. Applications}

Most one-dimensional quantum mechanical systems are described by
a Hamiltonian which is a densely defined unbounded self-adjoint
operator on $L^2(\Bbb R)$, having the form
$$
Hf(x)=-{1\over2}f^{\prime\prime}(x) + V(x)f(x)
$$
$V:\Bbb R\to\Bbb R$ being a continuous function representing the
potential.  If one wishes to simulate the behavior of such a quantum
system on a computer, one first has to discretize this differential
operator in an appropriate way.  Secondly, one has to develop
effective methods for computing with the discretized Hamiltonian.
For our purposes, we interpret the second goal to be that
of computing the spectrum.

In \cite{1}, \cite{2}, we discussed the problem of discretizing the above
Hamiltonians in such a way as to preserve the uncertainty principle.
We argued that the appropriate discretization is a bounded self-adjoint
operator of the form
$$
H_\sigma = -{1\over2}P^2_\sigma + V(Q_\sigma)\tag{5.1}
$$
where $P_\sigma$ and $Q_\sigma$ are the bounded self-adjoint operators
given by
$$
\align
P_\sigma f(x) &= {1\over{2i\sigma}}(f(x+\sigma)-f(x-\sigma)),\\
Q_\sigma f(x) &= {1\over\sigma}\sin(\sigma x)f(x).
\endalign
$$
Here, $\sigma$ is a small (positive rational) number representing the numerical
step size.

In even the simplest cases, one cannot carry out explicit calculations of the
spectrum of the discretized Hamiltonians 5.1
(see \cite{3},\cite{4},\cite{5},\cite{7},\cite{9}).  On the other hand, one can
apply the methods of the prceding sections in a straightforward manner to
obtain the spectrum of $H_\sigma$ as the limiting case of finite dimensional
eigenvalue distributions.

To see this, let $U$ and $V$ be the unitary operators
$$
\align
Uf(x)&=e^{i\sigma x}f(x),\\
Vf(x)&=f(x+2\sigma).
\endalign
$$
Then 5.1 can be written in the form $H_\sigma = \alpha A+\beta\I$, where
$\alpha$ and $\beta$ are real numbers and $A$ is a bounded self-adjoint
operator of the form
$$
A=V+V^* + v({1\over{2i}}(U-U^*)).\tag{5.2}
$$
Here, $v:\Bbb R\to\Bbb R$ is an appropriately rescaled version of $V$.
Thus, we are interested in the spectrum of $A$.  Of course, $A$ belongs
to the \cstar\ $C^*(U,V)$ generated by $U$ and $V$.  Following the
program outlined in section 4, we will first find a faithful representation
of $C^*(U,V)$ in which there is a compatible filtration, and then we will
look at the eigenvalues fo the approximating sequence of finite matrices.
Since $U$ and $V$ obey the commutation relation
$$
VU = e^{2i\sigma^2}UV\tag{5.3}
$$
and since $\sigma$ is a rational number, it follows that $C^*(U,V)$
is an irrational rotation \cstar.  Thus in order to define a faithful
representation, it is enough to specify a pair of unitary operators
$U_1$ and $V_1$ which satisfy the same commutation relation 5.3.  For
that, consider a Hilbert space spanned by an orthonormal basis
$\{e_n: n\in\Bbb Z\}$ and define $U_1$ and $V_1$ by
$$
\align
U_1e_n&=e^{-2in\sigma^2}e_n,\\
V_1e_n&=e_{n+1}.
\endalign
$$
The pair $\{U_1,V_1\}$ is irreducible, it obeys the required
commutation relations, and $A$ becomes the tridiagonal
operator $V_1 + V_1^* +D$, where $D$ is the diagonal operator
$$
De_n=v(-\sin(2\sigma^2 n))e_n,\qquad n\in\Bbb Z.
$$

Let $\Cal A = C^*(U_1,V_1)$.  Then the sequence of subspaces $H_1,H_2,\dots$
defined by $H_n=[e_{-n},e_{-n+1},\dots,e_n]$ is an $\Cal A$-filtration (see
example 2 of section 4).  Since concretely represented irrational rotation
\cstar s never contain nonzero compact operators, the spectrum of $A$
coincides with the essential spectrum of $A$.  Moreover, by Theorem 3.8,
the spectrum of $A$ is the set of essential points associated with the
sequence of matrices $\{A_n: n\geq1\}$, $A_n$ being the compression of $A$
to $H_n$.  The remaining points of $\Lambda$ are all transient points.
Finally,
since irrational rotation \cstar s are simple and have a unique trace, we may
conclude from Theorem 4.5 that for every $f\in C_0(\Bbb R)$ one has the
following convergence of the densities of eigenvalues,
$$
\lim_{n\to\infty}{1\over{2n+1}}(f(\lambda_1)+\dots +f(\lambda_{2n+1}) =
     \int_{-\infty}^{+\infty}f(x)\, d\mu_A(x),\tag{5.1}
$$
wehre $\mu_A$ is the spectral distribution of $A$ and
$\{\lambda_1,\dots,\lambda_{2n+1}\}$ is the eigenvalue list of $A_n$.

We remark that this representation is particularly convenient for
numerical eigenvalue computations using the QL or QR algorithms, because
the matrices $A_1,A_2,\dots$ are already in tridiagonal form relative to
the obvious basis.

\subheading{6.  Appendix. Failure of $\sigma_e(A) = \Lambda_e$}
In this section we give an example of a self-adjoint operator $A$ on
a Hilbert space spanned by an orthonormal basis $\{e_1,e_2,\dots\}$,
such that the relation of $A$ to the filtration
$H_n=[e_1,e_2,\dots,e_m]$, $n=1,2,\dots$
is pathological.  We show that if $A_n$ is the compression of
$A$ to $H_n$, then $\Lambda_e$  contains points not in the
spectrum of $A$, and in particular $\Lambda_e \neq \sigma_e(A)$.

Specifically, $A$ is a reflection (i.e., a self-adjoint unitary operator)
whose essential spectrum is $\{-1,+1\}$, whereas $0\in\Lambda_e$.  In the
example
$H = \l^2(\Bbb N)$, $e_n$ is the unit coordinate funtion at $n$, and
$A$ will be the operator in $\l^2(\Bbb N)$ induced by a permutation
$\pi$ of $\Bbb N$ satisfying $\pi^2 = \roman{id}$.  Let $I_n$ denote
the interval $\{1,2,\dots,n\}$, and let $\#S$ denote the number of
elements of a set $S$.  Suppose further that the permutation $\pi$ has
the following property
$$
\lim_{n\to\infty}\#(\pi(I_n)\setminus I_n) = +\infty.\tag{A.1}
$$
For each $n\geq1$, let $A_n$ be the compression of the operator
$A\xi(k) = \xi(\pi(k))$ to $H_n=\l^2(I_n)$.
Since $A$ maps $e_k$ to $e_{\pi(k)}$, it follows
that for every $k$ in the set $S_n = \{k\in I_n: \pi(k)\notin I_n\}$, we
have $A_ne_k=0$.
Since the cardinality of $S_n$ tends to infinity, it follows that
0 is an eigenvalue of every $A_n$ and that the multiplicity of this zero
eigenvalue tends to infinity.  Hence, 0 must belong to $\Lambda_e$.

Thus it suffices to exhibit an order two permutation having the property
A.1.
Let $E$ and $O$ denote respectively the set of even and odd integers in $\Bbb
N$.
In order to define a permutation $\pi$ with $\pi^2 = \roman{id}$, it is
enough to specify a bijection $f:E\to O$; once we have $f$ then we can define
$\pi$ by
$$
\pi(k) =
\cases
f(k), k\in E\\
f^{-1}(k), k\in O.
\endcases
$$
We will define a bijection $f:E\to O$ having the following property with
respect to the initial intervals $I_n$:
$$
\lim_{n\to\infty}\#(f(E\cap I_n)\setminus I_n) = +\infty.\tag{A.2}
$$
The resulting permutation will have the desired property A.1.

To construct $f$, split $E$ into two disjoint subsets $E_1\cup E_2$,
where $E_1 = 4\Bbb N$ is the set of multiples of 4 and $E_2$ is the
rest.  Both $E_1$ and $E_2$ are infinite of course.  First define \
$f$ on $E_1$ by
$$
f(k) = k^2 +1.\tag{A.3}
$$
Let $O_1=f(E_1)$, and put $O_2 = O\setminus O_1$.  Noting that $O_2$
is infinite (everything in $O_1$ is congruent to 1 mod 16, hence $O_2$
contains all the other odd numbers), we can complete the definition of
$f$ by choosing it to be any bijection of $E_2$ onto $O_2$.

To prove A.2, it suffices to show that $\#(f(E_1\cap I_n)\setminus I_n)$
tends to infinity with $n$.  Noting that
$$
E_1\cap I_n = \{4k : 1\leq k\leq [n/4]\},
$$
we see that
$$
f(E_1\cap I_n) = \{16k^2 + 1: 1\leq k\leq [n/4]\}.
$$
Thus every element in the subset $P_n$ of $E_1\cap I_n$ defined by
$$
P_n=\{4k: \frac{\sqrt n}{4} < k \leq \frac{n} {4}\}
$$
gets mapped by $f$ into the complement of $I_n$.  Hence
$$
\#(f(E_1\cap I_n)\setminus I_n)\geq \#P_n \geq {n\over4}(1- n^{-1/2}).
$$
This estimate implies that $\#(f(E\cap I_n)\setminus I_n) \to\infty$,
and A.2 follows.

\remark{Remark}  Notice that this misbehavior cannot be corrected by
averaging the eigenvalues of $A_n$ as we did in section 4.
Indeed, if $N_n(U)$ denotes the number of
eigenvalues of $A_n$ which belong to an open set $U$, then the above
estimate shows that for every neighborhood $U$ of 0 we have
$$
\liminf_{n\to\infty}{{N_n(U)}\over{\roman{dim}H_n}} \geq {1\over4}>0,
$$
in spite of the fact that 0 does not belong to the spectrum of $A$.
\endremark

%
%
%
%
%\newpage

\enddocument